\documentclass{LMCS}

\def\doi{9(2:3)2013}
\lmcsheading%
{\doi}
{1--14}
{}
{}
{Jun.~13, 2011}
{May.~22, 2013}
{}

\usepackage{enumerate}
\usepackage{hyperref}

\usepackage{amsmath,amssymb}
\usepackage{pgf}
\usepackage{tikz}
\theoremstyle{plain}

\begin{document}

\title[Infinite sequential Nash equilibrium]{Infinite sequential Nash equilibrium}

\author[S.~Le Roux]{St\'ephane Le Roux}	%required
\address{Mathematik AG1,
Technische Universit\"at Darmstadt
Schloßgartenstr. 7, D-64289 Darmstadt}	%required
\email{leroux@mathematik.tu-darmstadt.de}  %optional
%\thanks{thanks 1, optional.}	%optional

%% etc.

%% required for running head on odd and even pages, use suitable
%% abbreviations in case of long titles and many authors:

%% mandatory lists of keywords and classifications:
\ACMCCS{[{\bf  Theory of computation}]: Logic; [{\bf Mathematics of computing}]:  Continuous mathematics---Topology---Point-set topology; [{\bf  Applied computing}]: Operations research---Decision analysis;  Law, social and behavioral sciences---Economics} 

\begin{abstract}
  \noindent 
In game theory, the concept of Nash equilibrium reflects the collective stability of some individual strategies chosen by selfish agents. The concept pertains to different classes of games, \textit{e.g.} the sequential games, where the agents play in turn. Two existing results are relevant here: first, all finite such games have a Nash equilibrium (w.r.t. some given preferences) iff all the given preferences are acyclic; second, all infinite such games have a Nash equilibrium, if they involve two agents who compete for victory and if the actual plays making a given agent win (and the opponent lose) form a quasi-Borel set. This article generalises these two results \textit{via} a single result. More generally, under the axiomatic of Zermelo-Fraenkel plus the axiom of dependent choice (ZF+DC), it proves a transfer theorem for infinite sequential games: if all two-agent win-lose games that are built using a well-behaved class of sets have a Nash equilibrium, then all multi-agent multi-outcome games that are built using the same well-behaved class of sets have a Nash equilibrium, provided that the inverse relations of the agents' preferences are strictly well-founded. 
\end{abstract}

\maketitle

\section{Introduction}

Game theory is the theory of competitive interactions between decision makers having different interests. Its primary purpose is to further understand such real-world interactions through mathematical modelling. Apart from some earlier related works, the field of game theory is usually said to be born in the first part of the 20th century, especially thanks to von Neumann~\cite{NM44}, but also Borel~\cite{Borel21} and some others. Since then, it has been applied to many concrete areas such as economics, political science, evolutionary biology, computer science, \textit{etc}. Conversely, specific problems in these concrete areas have been triggering new general questions and have thus been helpful in developing game theory.

Surprisingly, game theory has also provided useful point of view and terminology to abstract areas such as logic and descriptive set theory. For instance, Martin~\cite{Martin90} relates (quasi-)Borel sets to the existence of winning strategies in some two-agent win-lose games, \textit{i.e.} games that involve only two players and two outcomes saying who wins. Conversely, one may wonder if and how such specific results may help develop general game theory. Especially, do such results scale up when considering the same game structure for many agents, many outcomes, and the usual generalisation of the notion of winning strategy, namely the notion of Nash equilibrium? 

This paper answers the question positively by proving a transfer theorem for infinite sequential games: if all two-agent win-lose games that are built using a well-behaved class of sets have winning strategies, then all (reasonable) multi-agent multi-outcome games that are built using the same well-behaved class of sets have Nash equilibria. Since the quasi-Borel sets are well-behaved and guarantee existence of winning strategies (as proved in~\cite{Martin90}), the application of the transfer theorem implies that all (reasonable) multi-agent multi-outcome games that are built using quasi-Borel sets have Nash equilibria. Martin's result is thus generalised by invoking Martin's result itself, which factors out most of the proof burden!

Section~\ref{subsect:not} sets some notations; Section~\ref{subsect:mr} describes the main result quickly; Section~\ref{subsect:pw} recalls related previous works; Section~\ref{subsect:prelim} recalls preliminary definitions and properties, mostly about well-foundedness; Section~\ref{sect:isg} introduces the main result in detail; Section~\ref{sect:iwl} applies the transfer theorem to quasi-Borel determinacy and to the axiom of determinacy.

\subsection{Notations}\label{subsect:not}

Let $C$ be a set. The following lists useful notations of well-known notions.

\begin{iteMize}{$\bullet$}
\item Outer universal quantifiers are sometimes omitted in formal statements.
\item In formal definitions, $x:=y$ means that $x$ is defined as $y$.
\item $C^*$ is the set of finite sequences over $C$, and $\epsilon$ is the empty sequence.
\item $|\gamma|$ is the length of $\gamma\in C^*$, and $\gamma=\gamma_0\gamma_1\dots\gamma_{|\gamma|-1}$.
\item $C^{2*}$ (resp. $C^{2*+1}$) is the set of finite sequences of even (resp. odd) length. 
\item $C^\omega$ is the set of infinite sequences over $C$. It is endowed with the product of the discrete topology on $C$.
\item For $\gamma,\gamma'\in C^*$, let $\gamma\cdot\gamma'$ (or $\gamma\gamma'$) denote the concatenation $\gamma_0\dots\gamma_{|\gamma|-1}\gamma'_0\dots\gamma'_{|\gamma'|-1}$.
\item For $\gamma\in C^*$ and $\gamma'\in C^\omega$, let $\gamma\gamma'$ denote the concatenation $\gamma_0\dots\gamma_{|\gamma|-1}\gamma'_0\gamma'_{1}\gamma'_{2}\dots$
\item For $X\subseteq C^*$ and $Y\subseteq C^*\cup C^\omega$, let $XY$ denote the concatenations of one element in $X$ and one element in $Y$.
\item $\gamma\sqsubseteq\gamma'$ means that $\gamma$ is a prefix of $\gamma'$ for $\gamma\in C^*$ and $\gamma'\in C^*\cup C^\omega$.
\item $C^{<n}:=\{\gamma\in C^*\,\mid\,|\gamma|<n\}$ are the sequences over $C$ of length shorter than $n$. 
\item $\gamma_{<n}$ is the prefix of length $n$ of $\gamma\in (C^*\cup C^\omega)\backslash C^{<n}$, in particular $\gamma_{<0}=\epsilon$.
\item $s:\subseteq C^*\to C$ means $s:dom(s)\to C$ where $dom(s)\subseteq C^*$.
\item $s\mid_{X}$ denotes the restriction of the function $s$ to $X$. 
\item If $s:\subseteq C^*\to C$ and $s':\subseteq C^*\to C$, then writing $s\cup s'$ requires that $dom(s)\cap dom(s')=\emptyset$, in which case $s\cup s':\subseteq C^*\to C$ is the function such that $s\cup s'(\gamma)$ is either not defined, or equal to $s(\gamma)$ if $\gamma\in dom(s)$, or equal to $s'(\gamma)$ if $\gamma\in dom(s')$.
\end{iteMize}

\subsection{The main result}\label{subsect:mr}

An infinite sequential game may involve infinitely many playing agents who play in turn. At each turn, one agent chooses an element in a given set of choices, and this specific choice determines the next agent who has to choose an element in the same set of choices, and so on. Each such an infinite sequence of choices is (objectively) summarised \textit{via} an outcome, and each agent gives the outcomes a (subjective) meaning by comparing them \textit{via} an individual binary relation.      

\begin{defi}[Infinite sequential game]\label{defn:ifg} An
  \emph{infinite sequential game}  is an object\\$\langle A,C,d,O,v,(\prec_a)_{a\in A}\rangle$ complying with the following.
\begin{iteMize}{$\bullet$}
\item $A$ is a non-empty set (of agents).
\item $C$ is a non-empty set (of choices).
\item $d:C^*\to A$ (assigns a decision maker to each stage of the game).
\item $O$ is a non-empty set (of possible outcomes of the game). 
\item $v:C^\omega\to O$ (uses outcomes to value the infinite sequences of choices).
\item Each $\prec_a$ is a binary relation over $O$ (modelling the preference of agent $a$).
\end{iteMize}
\end{defi}

Notice that letting the function $v$ from Definition~\ref{defn:ifg} map a clopen ball to the same outcome simulates a finite play ending with this outcome; also notice that a game with different sets of choices at each turn may be simulated by considering the union of these sets. Therefore, the above-defined infinite sequential games capture (most of) the competitive sequential processes that do not involve complex concepts such as time or knowledge.

In such a game, a strategy of an agent is a function saying what the agent would choose at every possible stage of the game where the agent is the actual decision maker. A strategy profile is a combination of one strategy per agent, as defined below.

\begin{defi}[Strategy]\label{defn:s}
A function $s:C^*\to C$ is called a \emph{strategy profile}. A function $s_a:d^{-1}(a)\to C$ is called a \emph{strategy for agent $a$}.
\end{defi}

Starting from the root of the game-tree and recursively following the unique choice made by the current decision maker yields an infinite sequence, or play, as defined below.

\begin{defi}[Induced play and outcome]\label{defn:ipo}
Let $s:C^*\to C$ be a strategy profile. The \emph{play $p=p(s)\in
  C^\omega$ induced by $s$} is defined inductively through its prefix:
$p_{n}:=s(p_{< n})$ for all $n\in\mathbb{N}$. One may also say that
$v\circ p(s)$ is the \emph{outcome induced by $s$}.
\end{defi}

The agents actually compare the strategy profiles \textit{via} their induced outcomes, like in Definition~\ref{defn:sne-ig} below. Notice that unlike the traditional definition of Nash equilibrium, the definition below uses two concepts: convertibility, \textit{i.e.} the ability of an agent to convert a strategy profile into another one, and preference. Furthermore, the preferences here are conceptually (slightly) different from the traditional notion of preference in game theory since they are not involved in the definition of Nash equilibrium in the same way. See~\cite{SLR09} for a detailed discussion on the matter.  

\begin{defi}[Sequential Nash equilibrium]\label{defn:sne-ig}
Let $\langle A,C,d,O,v,(\prec_a)_{a\in A}\rangle$ be a game. A
strategy profile $s:C^*\to C$ is a \emph{Nash equilibrium} if all
agents are stable, where an agent $a$ is \emph{stable} if he or she prefers
$v\circ p(s')$ over $v\circ p(s)$ for no $s'$ that coincides with $s$
outside of the subdomain $d^{-1}(a)$, said otherwise for no $s'$ into
which agent $a$ can convert $s$.
\[\begin{array}{r@{\quad:=\quad}l}
NE(s) & \forall a\in A,\quad stable(s,a)\\
stable(s,a) & \forall s'\in C^*\to C,\quad\neg(s\to_as'\,\wedge\,v\circ p(s)\prec_a v\circ p(s'))\\
s\to_a s' & s\mid_{C^*\setminus d^{-1}(a)}= s'\mid_{C^*\setminus d^{-1}(a)}
\end{array}\]
\end{defi}

Under the axiomatic of Zermelo-Fraenkel plus the axiom of dependent choice (ZF+DC), the theorem below is the main result of this article. (Note that Condition (4) in the theorem refers to two-agent games where one may either win or lose, and where determinacy means that one of the agents can win for sure, as will be defined formally in Sections~\ref{subsect:pw} and ~\ref{subsect:prelim}.)

\begin{thm}\label{thm:transfer-tree}
Let $\langle A,C,d,O,v,(\prec_a)_{a\in A}\rangle$ be an infinite sequential game. Let $\Gamma\subseteq\mathcal{P}(C^\omega)$ and assume the following.
\begin{enumerate}[\em(1)]
\item $O$ is well-orderable.
\item $\prec_a^{-1}$ is strictly well-founded for all $a\in A$.
\item $\forall O'\subseteq O,\forall (c,\gamma)\in C\times C^*,\quad v^{-1}(O')\cap\gamma(C-\{c\}) C^\omega\in\Gamma$ 
\item The game $\langle C,D,W\rangle$ is determined for all $W\in\Gamma$ and $D\subseteq C^*$.
\end{enumerate}
Then the game $\langle A,C,d,O,v,(\prec_a)_{a\in A}\rangle$ has a Nash equilibrium. 
\end{thm}

Conversely, it is foreseeable that a preference with an infinite
ascending chain may lead to games without equilibrium, as proved
later. Well-foundedness is therefore a key notion in game theory!

Conditions (3) and (4) above reflect the set-theoretic complexity of the distribution of the outcomes onto the infinite sequences of choices. For a large class of "simple" sets, \textit{i.e.} the quasi-Borel sets, a folklore extension of quasi-Borel determinacy will make Conditions (3) and (4) hold for countable sets $O$, thus generalising Martin's result~\cite{Martin90}.

\subsection{Previous works}\label{subsect:pw}

This section presents two existing results that are generalised by Theorem~\ref{thm:borel-nash}, as well a similar existing results. First of all, a Gale-Stewart game is an infinite sequential game involving just two agents who play alternately and two outcomes saying who wins. Let us call $a$ the agent who plays first, and $b$ the opponent.

\begin{defi}[Gale-Stewart games]
A \emph{Gale-Stewart game} is an object $\langle C,W\rangle$ where $C$ is a non-empty set (of possible choices at each step of the game) and $W$ is a subset of $C^\omega$ (and represents the winning set of agent $a$ among the infinite sequences over $C$, while $C^\omega\backslash W$ is the winning set of agent $b$).
\end{defi}

Following Definition~\ref{defn:s}, the strategies of agent $a$ (resp. $b$) are the functions of type $C^{2*}\to C$ (resp. $C^{2*+1}\to C$). The notion of induced play follows Definition~\ref{defn:ipo}.

In such specific games, it is natural to ask who wins, or said more carefully, whether one agent may actually win for sure.

\begin{defi}[Winning strategies and determined games]
Given a game $\langle C,W\rangle$, a strategy $s_a:C^{2*}\to C$ is a \emph{winning strategy for agent $a$} if $p(s_a\cup s_b)\in W$ for all $s_b:C^{2*+1}\to C$. Likewise, a strategy $s_b:C^{2*+1}\to C$ is a \emph{winning strategy for agent $b$} if $p(s_a\cup s_b)\notin W$ for all $s_a:C^{2*}\to C$. The game is \emph{determined} if either of the agents has a winning strategy.
\end{defi}

It turns out that the quasi-Borel sets, a large class of sets that are defined below (almost verbatim from \cite{Martin90}), enjoy nice determinacy properties such as in Theorem~\ref{thm:martin90}. 

\begin{defi}[Quasi-Borel]\label{defn:quasi-Borel}
Given a topological space $X$ the \emph{quasi-Borel sets} of $X$ form the smallest class $\Gamma$ of subsets of $X$ containing all open sets and closed under the operations:
\begin{enumerate}[(1)]
\item complementation;
\item countable union;
\item open-separated union.
\end{enumerate} 
Here $A$ comes from $\{B_j\,\mid\,j\in J\}$ by the operation of \emph{open-separated union}, or, equivalently $A$ is the \emph{open-separated union} of $\{B_j\,\mid\,j\in J\}$, if
\begin{enumerate}[a)]
\item $A=\cup\{B_j\,\mid\,j\in J\}$;
\item there are disjoint open sets $D_j$, $j\in J$, such that $B_j\subseteq D_j$ for each $j\in J$.
\end{enumerate}
\end{defi}

Notice that the two notions of Borel sets and quasi-Borel sets coincide for Polish spaces, \textit{e.g.} when $C$ is countable in Theorem~\ref{thm:martin90} below.

\begin{thm}[Martin~\cite{Martin90}]\label{thm:martin90}
Let $C$ be a non-empty set and $W$ be a quasi-Borel set of $C^\omega$. Then the game $\langle C,W\rangle$ is determined. 
\end{thm}

Actually, this paper does not only generalise quasi-Borel determinacy from descriptive set theory, but also a (much simpler) result~\cite{SLR-PhD08}\cite{SLR09} from multi-agent multi-outcome finite game theory: given agents, outcomes, and preferences from which one is allowed to build games, every finite sequential game has a Nash equilibrium if and only if the preferences of the agents are all acyclic. Note that this result is already a generalisation of existence of Nash equilibrium for real-valued finite games in extensive form, see~\cite{Kuhn53}, and existence of Nash equilibrium for finite games in extensive form with abstract outcomes and strict weak orders for agent preferences, see~\cite{OR94}. Likewise, quasi-Borel determinacy is already a generalisation of \cite{GS53}, \cite{Wolfe55}, and Borel determinacy~\cite{Martin75}.

Note that the idea of strengthening determinacy results for more general notions of equilibrium is not a novel idea: it has already been studied in parity games and Muller games, \textit{e.g.}, in \cite{Ummels05} and \cite{PS09}, respectively.

\subsection{Preliminaries}\label{subsect:prelim}

Let us first generalise slightly the notion of Gale-Stewart games, as is used in Theorem~\ref{thm:transfer-tree}, or equivalently Theorem~\ref{thm:wo-ne}. In these new games, in addition to the set of choices $C$ and a winning set $W$, the nodes where agent $a$ makes decisions are represented by a subset $D\subseteq C^*$ of finite sequences of choices. 

\begin{defi}[Infinite two-agent win-lose sequential games]\label{defn:i22s}
Such a game is an object $\langle C,D,W\rangle$ complying with the following.
\begin{iteMize}{$\bullet$}
\item $C$ is a non-empty set.
\item $D$ is a subset of $C^*$. ($D$ represents the stages of the game where agent $a$ makes decisions.)
\item $W$ is a subset of $C^\omega$.
\end{iteMize}
\end{defi}

Similarly as for Gale-Stewart games, strategies of agent $a$ are functions $s_a:D\to C$ and strategies of agent $b$ are functions $s_b:C^*\backslash D\to C$. The notion of induced play follows Definition~\ref{defn:ipo}, and winning strategies and determinacy are also defined similarly. Note that the Gale-Stewart games can be encoded as infinite two-agent win-lose sequential games (take $D:=C^{2*}$), which can in turn be encoded easily as infinite sequential games. Such an encoding is moreover faithful w.r.t. Nash equilibrium: in an infinite two-agent win-lose sequential game, $s_a\cup s_b$ is a Nash equilibrium iff either $s_a$ or $s_b$ is a winning strategy.

Let us now recall some set-theoretical results that are used in Section~\ref{sect:isg}. First, the following lemma is a well-known folklore result from set theory.

\begin{lem}[Well-ordered extension]\label{lem:woe}
Let $R\subseteq E\times E$ be a (strict) well-founded relation and assume that some $\prec_0$ well-orders $E$ (strictly). Then there exists a (strict) well-ordering $\prec$ of $E$ that includes $R$. 
\end{lem}

\proof
Following chapters 2D and 2G in~\cite{Moschovakis09}, let us define a rank function $\rho:E\to\lambda$ recursively, where $\lambda$ is an ordinal of cardinality $|R|$. So $\rho(x):=sup\{\rho(y)+1\,\mid\,yRx\}$. Let $\prec\subseteq E\times E$ be as follows.
\[x\prec y:=\rho(x)< \rho(y)\,\vee\,(\rho(x)=\rho(y)\,\wedge\,x\prec_0y)\] 
The relation $\prec$ includes $R$. Indeed, assume that $xRy$, so $\rho(x)+1\leq\rho(y)$ by definition of $\rho$, that is,  $\rho(x)<\rho(y)$. Therefore $x\prec y$ by definition of $\prec$.

The relation $\prec$ is linear order, as proved below. Assume that $x\prec y$ and $y\prec z$. If $\rho(x)<\rho(y)$ or $\rho(y)<\rho(z)$ then $\rho(x)<\rho(z)$, and therefore $x\prec z$ by definition. If $\rho(x)=\rho(y)=\rho(z)$ then $x\prec_0y$ and $y\prec_0z$, so $x\prec_0z$ by transitivity of $\prec_0$, so $x\prec z$ by definition. Moreover $\prec$ is irreflexive since $<$ and $\prec_0$ are irreflexive. Finally, if $x\neq y$ then either $\rho(x)<\rho(y)\,\vee\,\rho(y)<\rho(x)$ or $\rho(x)=\rho(y)$ (knowing $x\prec_0 y\,\vee\,y\prec_0 x$). Both cases yield $x\prec y\,\vee\,y\prec x$. So $\prec$ is a strict linear order. 

The relation $\prec$ is also well-founded: Let $S$ be a non-empty subset of $E$. Since ordinals are well-ordered by $<$, let $\alpha$ be the $<$-minimum of the set $\rho(S)$. So $\{x\in E\,\mid\,p(x)=\alpha\}\neq\emptyset$. Since $E$ is well-ordered by $\prec_0$, let $m$ be the $\prec_0$-minimum of $\{x\in E\,\mid\,p(x)=\alpha\}$. Now let $y\in S-\{m\}$. If $\alpha=\rho(m)<\rho(y)$ then $m\prec y$ by definition; if $\rho(m)=\rho(y)$, then $m\prec_0 y$ by construction, so $m\prec y$. Therefore $m$ is a $\prec$-minimum of $S$.
\qed

Terminal intervals (as opposed to initial) are defined below. They are sometimes also called upper sets in the literature.

\begin{defi}[Terminal interval]
Let $<$ be a strict linear order over $E$. The \emph{$<$-terminal intervals} are the subsets $I$ of $E$ that satisfy the formula $x<y\,\wedge\,x\in I\,\Rightarrow\,y\in I$. 
\end{defi}

\begin{obs}[Well-ordering of the terminal intervals]\label{obs:ti-wo}
Let $(E,<^{-1})$ be a strict well-order. Then the $<$-terminal intervals are (strictly) well-ordered by set-theoretic inclusion.  
\end{obs}

\section{Infinite sequential games}\label{sect:isg}

This section presents concepts and properties of increasing complexity, and culminates in the proof of the main theorem of the paper.

The notion of play induced by a strategy profile is extended below into the notion of possible plays under constraints, namely plays that are compatible with a partial strategy profile, as clarified by Lemma~\ref{lem:ppp}.\ref{lem:ppp6}. 

\begin{defi}[Possible plays under constraints]\label{defn:ppuc}
\[\forall s:\subseteq C^*\to C,\quad P(s):=\{p\in C^\omega\,\mid\,\forall n\in\mathbb{N},\,p_{<n}\in dom(s)\,\Rightarrow\,s(p_{<n})=p_{n}\}\]
\end{defi}

\begin{lem}\label{lem:ppp}
Let $C$ be a non-empty set.
\begin{enumerate}[\em(1)]
\item\label{lem:ppp0} $P(s\cup s')=P(s)\cap P(s')$, hence $s\subseteq s'$ implies $P(s')\subseteq P(s)$. 
\item\label{lem:ppp3} $\forall\gamma\in dom(s),\quad P(s\mid_{\gamma C^*})\cap\gamma C^\omega=P(s\mid_{\gamma s(\gamma)C^*})\cap\gamma s(\gamma)C^\omega$
\item\label{lem:ppp2} $\forall\gamma\notin dom(s),\quad P(s\mid_{\gamma C^*})\cap\gamma C^\omega=\bigcup_{c\in C}P(s\mid_{\gamma cC^*})\cap\gamma cC^\omega$
\item\label{lem:ppp4} $P(s\mid_{\gamma C^*})\cap\gamma C^\omega\subseteq P(s)$.
\item\label{lem:ppp6} $P(s)=\{p(t)\,\mid\,t:C^*\to C\,\wedge\, t\mid_{dom(s)}=s\}$
\end{enumerate} 
\end{lem}

\proof\hfill
\begin{enumerate}[(1)]
\item By definition.
\item Assume that $\gamma\in dom(s)$.
\[\begin{array}{c}
p\in P(s\mid_{\gamma C^*})\cap\gamma C^\omega\\
\Updownarrow\\
\gamma\sqsubseteq p\quad\wedge\quad\forall n\in\mathbb{N}, \gamma\sqsubseteq p_{<n}\,\wedge\,p_{<n}\in dom(s)\,\Rightarrow\,s(p_{<n})=p_{n}\\
\Updownarrow\\
\gamma\sqsubseteq p\quad\wedge\quad s(\gamma)=p_{|\gamma|}\quad\wedge\quad\forall n\in\mathbb{N}, \gamma s(\gamma)\sqsubseteq p_{<n}\,\wedge\,p_{<n}\in dom(s)\,\Rightarrow\,s(p_{<n})=p_{n}\\
\Updownarrow\\
\gamma s(\gamma)\sqsubseteq p\quad\wedge\forall n\in\mathbb{N}, \gamma s(\gamma)\sqsubseteq p_{<n}\,\wedge\,p_{<n}\in dom(s)\,\Rightarrow\,s(p_{<n})=p_{n}\\
\Updownarrow\\
p\in P(s\mid_{\gamma s(\gamma)C^*})\cap\gamma s(\gamma)C^\omega
\end{array}\]
\item Proof similar to the one above for Lemma~\ref{lem:ppp}.\ref{lem:ppp3}.
\item  By induction on $\gamma$ and Lemmas~\ref{lem:ppp}.\ref{lem:ppp2} and \ref{lem:ppp}.\ref{lem:ppp3}.
\item By double inclusion. Let $p\in P(s)$ and let us define $t:C^*\to C$ such that $t\mid_{dom(s)}:=s$ and $t(\gamma):=p_{|\gamma|}$ for all $\gamma\in C^*\backslash dom(s)$. By construction $t\mid_{dom(s)}:=s$ and $p(t)=p$. Conversely, let $t:C^*\to C$ such that $t\mid_{dom(s)}:=s$. Let $n\in\mathbb{N}$ and assume that $p(t)_{<n}\in dom(s)$, so $s(p(t)_{<n})=t(p(t)_{<n})=p(t)_n$ by assumption and by definition of the function $p$ respectively.\qedhere
\end{enumerate}

\noindent Given $a$ an agent, $\gamma$ a finite sequence of choices, and $s$ a strategy profile, let us consider the plays starting with prefix $\gamma$ and satisfying the constraint given by $s\mid_{d^{-1}(a)}$ in the subgame rooted at node $\gamma$ (so that there is at least one such play). This subset of plays is pointwise mapped to a subset of outcomes by the valuation function $v$. The smallest "upper set" (according to the preference of $a$) including this subset of outcomes represents the guarantee of the agent $a$ under the above constraints. This is formally defined below.

\begin{defi}[Agent guarantee]\label{defn:ag}
 Let $\langle A,C,d,O,v,(<_a)_{a\in A}\rangle$ be a game where the $<_a$ are strict linear orders, and where $x\leq_ay$ means $x<_ay\vee x=y$.
\[\begin{array}{l}
\forall a\in A,\forall \gamma\in C^*,\forall s: C^*\to C,\\
g_a(\gamma,s):=\{o\in O\,\mid\,\exists p\in P(s\mid_{d^{-1}(a)\cap\gamma C^*})\cap\gamma C^\omega,\, v(p)\leq_a o\}
\end{array}\]
\end{defi}

\begin{lem}\label{lem:g}
Let $\langle A,C,d,O,v,(<_a)_{a\in A}\rangle$ be a game where the $<_a$ are strict linear orders.
\begin{enumerate}[\em(1)]
\item\label{lem:g1} $g_a(\gamma,s)$ is the smallest $<_a$-terminal interval including $v(P(s\mid_{d^{-1}(a)\cap\gamma C^*})\cap\gamma C^\omega)$. 
\item\label{lem:g3} If $d(\gamma)= a$ then $g_a(\gamma,s)=g_a(\gamma s(\gamma),s)$.
\item\label{lem:g2} If $d(\gamma)\neq a$ then $g_a(\gamma,s)=\cup_{c\in C}g_a(\gamma c,s)$.
\item\label{lem:g5} If $s\mid_{C^*\backslash\gamma\gamma'C^*}=t\mid_{C^*\backslash\gamma\gamma'C^*}$ and $g_a(\gamma\gamma',s)\subseteq g_a(\gamma\gamma',t)$, then $g_a(\gamma,s)\subseteq g_a(\gamma,t)$.
\end{enumerate}
\end{lem}

\proof\hfill
\begin{enumerate}[(1)]
\item Let $(o,o')\in g_a(\gamma,s)\times O$ be such that $o<_a o'$ and let $I$ be a $<_a$-terminal interval including $v(P(s\mid_{d^{-1}(a)\cap\gamma C^*})\cap\gamma C^\omega)$. Let $p\in P(s\mid_{d^{-1}(a)\cap\gamma C^*})\cap\gamma C^*$ such that $v(p)\leq_a o$, by Definition~\ref{defn:ag}. On the one hand, $v(p)\leq_a o'$ by transitivity, so $o'\in g_a(\gamma,s)$, which shows that $g_a(\gamma,s)$ is a $<_a$-terminal interval; on the other hand, since $v(p)\in I$ by inclusion assumption, $o$ is also in $I$ by terminal assumption.
\item By rewriting Definition~\ref{defn:ag} using Lemma~\ref{lem:ppp}.\ref{lem:ppp3}.
\item By rewriting Definition~\ref{defn:ag} using Lemma~\ref{lem:ppp}.\ref{lem:ppp2}.
\item By induction on $\gamma'$. The base case $\gamma'=\epsilon$ is clear. For the induction step, assume that $s\mid_{C^*\backslash\gamma\gamma' cC^*}=t\mid_{C^*\backslash\gamma\gamma' cC^*}$ for some $c\in C$ and that $g_a(\gamma\gamma'c,s)\subseteq g_a(\gamma\gamma'c,t)$. Then $g_a(\gamma\gamma',s)\subseteq g_a(\gamma\gamma',t)$ by case-splitting on whether $d(\gamma\gamma')=a$ (and in case yes, whether $c=s(\gamma\gamma')$), by both assumptions, and by Lemmas~\ref{lem:g}.\ref{lem:g2} and \ref{lem:g}.\ref{lem:g3}. Therefore $g_a(\gamma,s)\subseteq g_a(\gamma,t)$ by induction hypothesis.\qed\smallskip
\end{enumerate}

\noindent Assume a finite sequence of choices $\gamma$ as fixed constraint. By choosing different strategies as additional constraint, the agent $a$ may obtain different guarantees. Informally, since the agent may prefer smaller upper sets of guarantee, the infimum of these available guarantees, \textit{i.e.} their intersection, is a relevant notion defined below. 

\begin{defi}[Best guarantee]\label{defn:bg}
Let $\langle A,C,d,O,v,(<_a)_{a\in A}\rangle$ be a game where the $<_a$ are strict linear orders.
\[\forall a\in A,\forall \gamma\in C^*,\quad G_a(\gamma):=\bigcap_{s:C^*\to C}g_a(\gamma,s)\]
\end{defi}

The lemma below says that, under a well-ordering condition, an agent can maximise his or her guarantee, \textit{i.e.} secure the best possible guarantee. 

\begin{lem}\label{lem:bgf}
Let $\langle A,C,d,O,v,(<_a)_{a\in A}\rangle$ be a game. If the $<_a^{-1}$ are strict well-orders, then the following holds.
\[\forall a\in A,\forall \gamma\in C^*,\exists\mu_{a}^{\gamma}: C^*\to C,\quad g_a(\gamma,\mu_{a}^{\gamma})= G_a(\gamma)\]
\end{lem}

\proof
Let $a\in A$ and $\gamma\in C^*$. By Lemma~\ref{lem:g}.\ref{lem:g1} the $g_a(\gamma,s)$ are non-empty terminal intervals for all $s$. By Observation~\ref{obs:ti-wo} these terminal intervals are well-ordered by inclusion, so the (non-empty) set $\{g_a(\gamma,s)\,\mid\,s:C^*\to C\}$ has a minimum $g_a(\gamma,\mu_{a}^{\gamma})$ for some witness $\mu_{a}^{\gamma}$, and the minimum equals $G_a(\gamma)$ by Definition~\ref{defn:bg}.
\qed

Given an infinite sequential game, the definition below builds a special strategy profile as follows. Start with a dummy strategy profile and let the root owner change his or her own strategy to maximise his or her guarantee at the root, unless it is already maximal. Now move to the node of the infinite tree that is (newly) chosen by the root owner, and let the new node owner maximise his or her guarantee at this node, and so on.

\begin{defi}[Deepening the guarantee]\label{defn:dg}
Let $\langle A,C,d,O,v,(<_a)_{a\in A}\rangle$ be a game where the
$<_a^{-1}$ are strict well-orders. Let $s_0:C^*\to C$ be,
\textit{e.g.,} a constant strategy profile, and let us define a
sequence of strategy profiles inductively: let us assume that $s_n$ is
defined, let $\gamma:=p(s_n)_{<n}$, let $a:=d(\gamma)$, and let us
define $s_{n+1}$ as follows.
\begin{iteMize}{$\bullet$}
\item If $g_a(\gamma,s_n)=G_a(\gamma)$ then $s_{n+1}:=s_n$.
\item Otherwise $s_{n+1}:=s_n\mid_{C^*\backslash (d^{-1}(a)\cap \gamma C^*)}\cup \mu_{a}^{\gamma}\mid_{d^{-1}(a)\cap \gamma C^*}$\\
where $\mu_{a}^{\gamma}$ is a witness from Lemma~\ref{lem:bgf}.
\end{iteMize}
\end{defi}

As proved below, the deepening process converges to a strategy profile that maximises guarantees along its play, at each node and for the node owner.  Note that the process modifies strategies only when necessary, otherwise it could fail even for a one-player win-lose game: consider the game $\langle \{0,1\},\{0,1\}^*,\{0,1\}^\omega-\{0^\omega\}\rangle$, which matches Definition~\ref{defn:i22s}. It is basically played by one single agent who wins if managing not to choose the zero sequence. If the deepening process allowed modifications of optimal strategies, for all $n$ the agent could choose some $s_n$ inducing $0^n1^\omega$, and the limit strategy would make the agent lose. 

\begin{lem}[Properties of the deepening]\label{lem:pd}
Recall $(s_n)_{n\in\mathbb{N}}$ from Definition~\ref{defn:dg}.
\begin{enumerate}[\em(1)]
\item\label{lem:pd1} $s_{n+k}\mid_{C^*\backslash p(s_n)_{<n}C^*}=s_{n}\mid_{C^*\backslash p(s_n)_{<n}C^*}$
\item\label{lem:pd2} Let $\sigma:C^*\to C$ be such that $\sigma(\gamma):=s_{|\gamma|+1}(\gamma)$ for all $\gamma\in C^*$. The sequence $(s_{n})_{n\in\mathbb{N}}$ converges towards $\sigma$; more specifically $s_{n+k}\mid_{C^*\backslash p(s_n)_{<n}C^*}=\sigma\mid_{C^*\backslash p(s_n)_{<n}C^*}$. 
\item\label{lem:pd3} $p(s_{n+k})_{<n}=p_{<n}$, where $p:=p(\sigma)$.\\
(This result and notation are often used implicitly in the rest of the article.)
\item\label{lem:pd-1} $d(p_{<n})=a\quad\Rightarrow\quad g_a(p_{<n},s_{n+1})=G_a(p_{<n})$ 
\item\label{lem:pd4} $\forall a\in A,\quad g_a(p_{<n+1},s_{n+1})\subseteq g_a(p_{<n},s_{n+1})\subseteq g_a(p_{<n},s_n)$
\item\label{lem:pd5} $d(p_{<n})=a\quad\Rightarrow\quad g_a(p_{<n},s_{n+1+k})= G_a(p_{<n})$.
\item\label{lem:pd6} $d(p_{<n})=a\quad\Rightarrow\quad g_a(p_{<n},\sigma)= G_a(p_{<n})$
\end{enumerate}
\end{lem}

\proof\hfill
\begin{enumerate}[(1)]
\item First prove $p(s_{n+k})_{<n}=p(s_n)_{<n}$ by induction on $k$ and Definition~\ref{defn:ipo}, then use it to prove the claim by induction on $k$.
\item By Lemma~\ref{lem:pd}.\ref{lem:pd1}.
\item By  Lemma~\ref{lem:pd}.\ref{lem:pd2}.
\item By case splitting along Definition~\ref{defn:dg} and using Definition~\ref{defn:ag}.
\item If $d(p_{<n})=a$, let $h:=g_a(p_{<n+1},s_{n+1})$; so $h=g_a(p(s_{n+1})_{<n+1},s_{n+1})$ by Lemma~\ref{lem:pd}.\ref{lem:pd3}, so $h=g_a(p(s_{n+1})_{<n}s_{n+1}(p_{<n}),s_{n+1})$ by Definition~\ref{defn:ipo} of the induced play, so $h=g_a(p(s_{n+1})_{<n},s_{n+1})$ by Lemma~\ref{lem:g}.\ref{lem:g3}, so $h=g_a(p_{<n},s_{n+1})$ by Lemma~\ref{lem:pd}.\ref{lem:pd3} again, so $h=G_a(p_{<n})$ by Lemma~\ref{lem:pd}.\ref{lem:pd-1}, so $h\subseteq g_a(p_{<n},s_n)$ by Definition~\ref{defn:bg}. Alternately, if $d(p_{<n})\neq a$, then $g_a(p_{<n+1},s_{n+1})\subseteq g_a(p_{<n},s_{n+1})=g_a(p_{<n},s_n)$ respectively by Lemma~\ref{lem:g}.\ref{lem:g2} and since $s_{n+1}\mid_{d^{-1}(a)}=s_n\mid_{d^{-1}(a)}$ by Definition~\ref{defn:dg}.
\item By induction on the number of $j$ between (arbitrary) $n+1$ and $n+k$ such that $s_{j}\mid_{d^{-1}(a)}\neq s_{j+1}\mid_{d^{-1}(a)}$, \textit{i.e.} intuitively, the number of times that agent $a$ refines its strategy. Base step, there is no such $j$, so $s_{n+1}\mid_{d^{-1}(a)}=s_{n+1+k}\mid_{d^{-1}(a)}$, so $g_a(p_{<n},s_{n+1+k})=g_a(p_{<n},s_{n+1})= G_a(p_{<n})$ by Definition~\ref{defn:ag} and Lemma~\ref{lem:pd}.\ref{lem:pd-1}, respectively. Inductive step, let $j$ be the minimal such $j$, so $g_a(p_{<j}, s_{n+1+k})=G_a(p_{<j})\subseteq g_a(p_{<j}, s_{j})$ by induction hypothesis and Definition~\ref{defn:bg}, respectively. It implies $g_a(p_{<n}, s_{n+1+k})\subseteq g_a(p_{<n}, s_{j})$ by Lemma~\ref{lem:g}.\ref{lem:g5}, since $s_{n+1+k}\mid_{C^*\backslash p_{<j}C^*}=s_{j}\mid_{C^*\backslash p_{<j}C^*}$ by Lemma~\ref{lem:pd}.\ref{lem:pd1}. But $g_a(p_{<n}, s_{j})=g_a(p_{<n}, s_{n+1})=G(p_{<n})$ respectively since $s_{j}\mid_{d^{-1}(a)}=s_{n+1}\mid_{d^{-1}(a)}$ (by definition of $j$) and by Lemma~\ref{lem:pd}.\ref{lem:pd-1}, so $g_a(p_{<n}, s_{n+1+k})\subseteq G(p_{<n})$.
\item By Lemmas~\ref{lem:pd}.\ref{lem:pd4} and~\ref{lem:g}.\ref{lem:g1} and Observation~\ref{obs:ti-wo}, the sequence $(g_a(p_{<k},s_k))_{k\in\mathbb{N}}$ is constant from some $k_0$ onwards. (Intuitively, no agent can keep on refining his/her strategy infinitely often -when moving deeper and deeper along the play $p$ being defined- because each refinement involves a strict improvement according to a preference that has no infinite ascending chain.) Because of Definition~\ref{defn:dg}, it implies that the sequence $(s_k\mid_{d^{-1}(a)})$ is constant from $k_0$ onwards, therefore equal to $\sigma\mid_{d^{-1}(a)}$. So $g_a(p_{<n},\sigma)=g_a(p_{<n},s_{n+k_0+1})$, which is equal to $G_a(p_{<n})$ by Lemma~\ref{lem:pd}.\ref{lem:pd5}.\qed\smallskip
\end{enumerate}

\noindent The proof of the theorem below builds a Nash equilibrium that shares the same induced play $p$ with the limit strategy profile $\sigma$ from Lemma~\ref{lem:pd}. The idea is that at each node along $p$, the node owner is threatened by the opponents collectively that deviating from $p$ would not yield any better outcome than $v(p)$. Such local, collective, threatening strategies exist for two combined reasons: by Lemma~\ref{lem:pd}.\ref{lem:pd6} at a given node along $p$, the node owner cannot secure any better outcome than $v(p)$, so by determinacy assumption, the opponents can collectively exclude the outcomes that are better than $v(p)$ for the node owner.

\begin{thm}\label{thm:wo-ne}
Let $\langle A,C,d,O,v,(\prec_a)_{a\in A}\rangle$ be an infinite sequential game. Let $\Gamma\subseteq\mathcal{P}(C^\omega)$ and assume the following.
\begin{enumerate}[\em(1)]
\item\label{hyp:wo-ne-woa} $O$ is well-orderable.
\item\label{hyp:wo-ne-wf} $\prec_a^{-1}$ is strictly well-founded for all $a\in A$.
\item\label{hyp:wo-ne-stable} $\forall O'\subseteq O,\forall (c,\gamma)\in C\times C^*,\quad v^{-1}(O')\cap\gamma(C-\{c\}) C^\omega\in\Gamma$ 
\item The game $\langle C,D,W\rangle$ is determined for all $W\in\Gamma$ and $D\subseteq C^*$.
\end{enumerate}
Then the game $\langle A,C,d,O,v,(\prec_a)_{a\in A}\rangle$ has a Nash equilibrium. 
\end{thm}

\proof
It suffices to prove the claim when the $\prec_a^{-1}$ are strict well-orders. Indeed, by Lemma~\ref{lem:woe} and Assumptions~\ref{hyp:wo-ne-woa} and \ref{hyp:wo-ne-wf}, each of the $\prec_a^{-1}$ is included in some strict well-order, and any Nash equilibrium using the strict well-orders as preferences would be also a Nash equilibrium according to the original preferences, by Definition~\ref{defn:sne-ig}. So assume that the $\prec_a^{-1}$ are strict well-orders.

Recall $\sigma$ and $p=p(\sigma)$ from Lemma~\ref{lem:pd}. At each node of the play $p$, let us define one win-lose game involving the owner of the node versus the other agents. Let $n\in\mathbb{N}$, let $a:=d(p_{<n})$ and let $g_n:=\langle C,D_n,W_n\rangle$ where $D_n:=d^{-1}(a)\cup (C^*\backslash p_{<n}C^*)$ and $W_n:=\{p'\in p_{<n}(C-\{p_n\})C^\omega\,\mid\,v(p)\prec_a v(p')\}$. Since $W_n=p_{<n}(C-\{p_n\})C^\omega\cap v^{-1}\{o\in O\,\mid\,v(p)\prec_a o\}$, it is in $\Gamma$ by Assumption~\ref{hyp:wo-ne-stable}. So the game $g_n$ is determined. 

Now let us prove that agent $a$ loses the game $g_n$. By contradiction, let $s$ be a strategy profile such that $s\mid_{D_n}$ is a winning strategy for agent $a$ in the game $g_n$. Let $h:=P(s\mid_{d^{-1}(a)\cap p_{<n}C^*})\cap p_{<n}C^\omega$, so that $h=P(s\mid_{D_n\cap p_{<n}C^*})\cap p_{<n}C^\omega$ since $d^{-1}(a)\cap p_{<n}C^*=D_n\cap p_{<n}C^*$ by definition of $D_n$, so $h\subseteq P(s\mid_{D_n})$ by Lemma~\ref{lem:ppp}.\ref{lem:ppp4}. Since $s\mid_{D_n}$ is winning, $P(s\mid_{D_n})\subseteq W_n$ by Lemma~\ref{lem:ppp}.\ref{lem:ppp6}, and $W_n\subseteq v^{-1}\{o\in O\,\mid\,v(p)\prec_a o\}$ by definition of $W_n$. Combining these equalities and inclusions yields $P(s\mid_{d^{-1}(a)\cap p_{<n}C^*})\cap p_{<n}C^\omega\subseteq v^{-1}\{o\in O\,\mid\,v(p)\prec_a o\}$, so $G_a(p_{<n})\subseteq g_a(p_{<n},s)\subseteq \{o\in O\,\mid\,v(p)\prec_a o\}$ by Definition~\ref{defn:bg} and by Lemma~\ref{lem:g}.\ref{lem:g1}, respectively. Now recall that $p:=p(\sigma)$, so $p\in P(\sigma\mid_{d^{-1}(a)\cap p_{<n}C^*})$ by Definition~\ref{defn:ppuc}, so $v(p)\in g_a(p_{<n},\sigma)= G_a(p_{<n})$ respectively by Definition~\ref{defn:ag} and Lemma~\ref{lem:pd}.\ref{lem:pd6}. Since $v(p)\notin\{o\in O\,\mid\,v(p)\prec_a o\}$ by definition, it contradicts the above inclusion. So let $t_n$ be a winning strategy for the opponent of agent $a$ in the game $g_n$, so that $P(t_n)\cap W_n=\emptyset$ by Lemma~\ref{lem:ppp}.\ref{lem:ppp6}.

Now let us construct a strategy profile by appending the local threats $t_n$ to the special play $p$. Let us define $s$ as follows, and note that $p(s)=p$.

\begin{iteMize}{$\bullet$}
\item $\forall n\in\mathbb{N},\quad s(p_{<n}):=p_{n}$
\item $\forall n\in\mathbb{N},\forall c\in C-\{p_{n}\},\forall\gamma\in C^*,\quad s(p_{<n}c\gamma):=t_n(p_{<n}c\gamma)$ if $d(p_{<n}c\gamma)\neq d(p_{<n})$ and $s(p_{<n}c\gamma):=c$ otherwise
\end{iteMize}

\noindent Let us eventually prove that $s$ is a Nash equilibrium, according to Definition~\ref{defn:sne-ig}. Let $a$ be an agent, let $s':C^*\to C$ be such that $s'\mid_{C^*\backslash d^{-1}(a)}=s\mid_{C^*\backslash d^{-1}(a)}$, let $p':=p(s')$, and let us show that $\neg(v(p)\prec_a v(p'))$. If $p'=p$ then $\neg (v(p)\prec_a v(p'))$ by irreflexivity of strict well-orders, so now assume that $p'\neq p$. Let $m$ be such that $p'_{<m}=p_{<m}$ and $p'_m\neq p_m$, so $s'(p_{<m})\neq s(p_{<m})$ by Definition~\ref{defn:ipo}. Since $s'\mid_{C^*\backslash d^{-1}(a)}=s\mid_{C^*\backslash d^{-1}(a)}$ by assumption, we have $d(p_{<m})=a$, so $s'\mid_{p'_{<m+1}C^*\backslash d^{-1}(a)}=s\mid_{p'_{<m+1}C^*\backslash d^{-1}(a)}=t_m\mid_{p'_{<m+1}C^*}$ by definition of $s$, so $\{p'\}=P(s')\subseteq P(t_m)$ by Lemmas~\ref{lem:ppp}.\ref{lem:ppp6} and \ref{lem:ppp}.\ref{lem:ppp0}. Moreover recall that $P(t_m)\cap W_m=\emptyset$, so $p'\in p_{<m}(C-\{p_m\})C^\omega\backslash W_m$ by definition of $m$, so $v(p')\prec_a v(p)$ by definition of $W_m$.  
\qed

\section{Applications of the transfer theorem}\label{sect:iwl}

The two applications in this section, Theorems~\ref{thm:borel-nash} and \ref{thm:det-nash}, are two "extreme" corollaries of Theorem~\ref{thm:wo-ne}. First, Theorem~\ref{thm:borel-nash} generalises the quasi-Borel determinacy~\cite{Martin90} (which is proved in ZFC). Second, Theorem~\ref{thm:det-nash} assumes the determinacy of all sets (which is inconsistent with ZFC) and derives existence of Nash equilibrium for all infinite sequential games, provided that the inverse relations of the agents' preferences are strictly well-founded. Between these two extremes, one may also instantiate Theorem~\ref{thm:wo-ne} with classes $\Gamma$ of various complexity together with determinacy assumptions for $\Gamma$. 

Before stating Theorem~\ref{thm:borel-nash}, let us address a slight issue, and since the solution of the issue belongs to the folklore, the related definitions and proofs will be only sketched. The issue is: Martin~\cite{Martin90} formally proved determinacy of quasi-Borel sets only for games where agents play strictly alternately, \textit{i.e.} in games $\langle C,W\rangle$, whereas Theorem~\ref{thm:wo-ne} requires a determinacy statement for games $\langle C,D,W\rangle$.

The folklore solution is: For each game where the two agents play in no specific order, let us translate it in a new game where the agents play alternately as follows. Each time the same agent would play twice in a row, let us force the opponent to play a dummy move in between. This dummy insertion transforms each finite/infinite sequence of moves in the original game into a new finite/infinite sequence of moves. The new winning set of an agent is the image of his/her old winning set by the insertion function, union a set that makes the opponent lose when refusing to play the dummy moves as prescribed. Strategies are naturally translated along the insertion, and a winning strategy in the new (Gale-Stewart) game is eventually translated back into a winning strategy in the original game. 

More formally, let $g:=\langle C,D,W\rangle$ be a game, assume without loss of generality that $\epsilon\in D$, and let us define an insertion function $\iota:C^*\to C^*$ inductively as follows, where $c_0\in C$ is fixed and $\gamma\in C^*$ and $c\in C$ are arbitrary.

\begin{iteMize}{$\bullet$}
\item $\iota(\epsilon):=\epsilon\quad$
\item $\iota(\gamma c):=\iota(\gamma)cc_0\quad$ if $\gamma\in D\Leftrightarrow\gamma c\in D\quad$ 
\item $\iota(\gamma c):=\iota(\gamma)c\quad$ if $\gamma\in D\not\Leftrightarrow\gamma c\in D\quad$
\end{iteMize}

Let us extend $\iota$ to infinite sequences by $\iota(\alpha)_n:=\iota(\alpha_{< n})_n$ for all $\alpha\in C^\omega$ and $n\in\mathbb{N}$. Many useful properties of $\iota$ carry over from finite to infinite sequences: $\iota$ is injective and computable. Also, $\iota(C^\omega)$ is closed and $\iota(\gamma C^\omega)=\iota(\gamma)C^\omega\cap\iota(C^\omega)$ for all $\gamma\in C^*$.

After "stretching" the game $g$ by using $\iota$ to insert dummy nodes, let us add the fatal nodes $F:=\{\iota(\gamma)cd\,\mid\,\iota(\gamma c)=\iota(\gamma)cc_0\,\wedge\,d\in C-\{c_0\}\}$. These are the nodes that an agent would hit just after not playing as prescribed at a dummy node.

Let us now split $C^\omega\backslash\iota(C^\omega)$ into two disjoint open sets $W_F$ and $L_F$ as below, respectively making the agents $b$ and $a$ lose when not playing the dummy moves as prescribed.
\[W_F:=\bigcup_{\gamma\in F\cap C^{2*}}\gamma C^\omega\qquad L_F:=\bigcup_{\gamma\in F\cap C^{2*+1}}\gamma C^\omega\]
Note that every winning strategy for a player in $g':=\langle C,\iota(W)\cup W_F\rangle$ can be straightforwardly converted into a winning strategy for the same player in $g$. Finally, to show that $g'$ is determined, it would suffice to show that $\iota(W)\cup W_F$ is quasi-Borel whenever $W$ is quasi-Borel. This last step might actually be also doable for other interesting classes of sets, hence the following abstracted summary, which is then invoked in the proof of Theorem~\ref{thm:borel-nash}.

\begin{lem}\label{lem:no-gsg}
Let $C\neq\emptyset$, let $\Gamma\subseteq\mathcal{P}(C^\omega)$, and assume the following. 
\begin{iteMize}{$\bullet$}
\item If $S\in\Gamma$ and $E$ is an open set of $C^\omega$, then $\iota(S)\cup E\in\Gamma$. 
\item The game $\langle C,W\rangle$ is determined for all $W\in\Gamma$.
\end{iteMize}
Then the game $\langle C,D,W'\rangle$ is also determined for all $W'\in\Gamma$ and $D\subseteq C^*$.
\end{lem}

\begin{thm}\label{thm:borel-nash}
Let $A$ be non-empty, let $C$ have at least two elements, let $O$ be non-empty countable, and for all $a\in A$, let $\prec_a$ be a binary relation over $O$. Then the following two propositions are equivalent.
\begin{enumerate}[\em(1)]
\item $\prec_a^{-1}$ is strictly well-founded for all $a\in A$.
\item For all $d:C^*\to A$ and $v:C^\omega\to O$, if $v^{-1}(o)$ is quasi-Borel for all $o\in O$, the game $\langle A,C,d,O,v,(\prec_a)_{a\in A}\rangle$ has a Nash equilibrium. 
\end{enumerate}
\end{thm}

\proof
Let us first prove $2\Rightarrow 1$ by contraposition. Assume that $\prec_a^{-1}$ is not strictly well-founded for some $a\in A$, so $\prec_a$ has an infinite ascending chain $(o_n)_{n\in\mathbb{N}}$, that is, $o_n\prec_a o_{n+1}$ for all $n\in\mathbb{N}$. Let $d(\gamma):=a$ for all $\gamma\in C^*$. Let $c\in C$, let $v(c^\omega):=o_0$ and let $v(c^n\alpha):=o_{n+1}$ for all $n\in\mathbb{N}$ and for all $\alpha\in (C-\{c\})C^\omega$ (infinite sequences that do not start with the choice $c$). For all $o\in O$, $v^{-1}(o)$ is quasi-Borel since it is either empty, or an open ball, or the closed set $\{c^\omega\}$. The game $\langle A,C,d,O,v,\prec_a\rangle$ has no Nash equilibrium nonetheless, since every strategy (profile) of agent $a$ induces one $o_k$ for some $k$ and may be improved upon by a strategy profile inducing $o_{k+1}$. 

$1\Rightarrow 2$ may be proved by Theorem~\ref{thm:wo-ne} where $\Gamma$ is instantiated with the quasi-Borel sets, so let us check whether the hypotheses hold. First hypothesis, $O$ is countable and therefore well-orderable (by the enumeration that witnesses its countability). Second hypothesis, by assumption $\prec_a^{-1}$ is strictly well-founded for all $a\in A$. Third hypothesis, let $O'\subseteq O$ and $\gamma\in C^*$. Since $v^{-1}(O')=\bigcup\{v^{-1}(o)\,\mid\,o\in O'\}$ is a countable union of sets that are quasi-Borel by assumption, it is also quasi-Borel by definition, so by intersection with an open ball, $v^{-1}(O')\cap\gamma C^\omega$ is still quasi-Borel by definition. The fourth hypothesis may be proved by Lemma~\ref{lem:no-gsg}, so let us check whether the hypotheses hold. The second hypothesis was proved by Martin~\cite{Martin90}, so it suffices to prove that the quasi-Borel sets are closed under the insertion function $\iota$. Let us proceed by transfinite induction, which is possible thanks to the quasi-Borel hierarchy that is described just before Lemma 1.1 in~\cite{Martin90}.
\begin{iteMize}{$\bullet$}
\item Base case, let $W$ be an open set, so $W$ is the union $\cup_{j\in J} \gamma_jC^\omega$ of some basic clopen sets. Since $\iota(\gamma_jC^\omega)=\iota(\gamma_j)C^\omega\cap\iota(C^\omega)$ by property of $\iota$, the set $\iota(W)=\cup_{j\in J}\iota(\gamma_jC^\omega)=(\cup_{j\in J}\iota(\gamma_j)C^\omega)\cap\iota(C^\omega)$ is the intersection of an open set and a closed set, so it is quasi-Borel.
\item Let $W$ come from complementation of $C^\omega\backslash W$, so by induction hypothesis $\iota(C^\omega\backslash W)$ is quasi-Borel. Since $\iota(C^\omega)$ is closed and $\iota(W) = \iota(C^\omega)\backslash\iota(C^\omega\backslash W)$ by injectivity, it is quasi-Borel.
\item Let $W$ come from a countable union $\cup_{n\in\mathbb{N}}W_n$, so $\iota(W)=\cup_{n\in\mathbb{N}}\iota(W_n)$ is a countable union of sets that are quasi-Borel by induction hypothesis, so it is quasi-Borel, too.
\item Let $W$ come from an open-separated union $\cup_{j\in J}W_j$ that is witnessed by $W_j\subseteq\cup_{i\in I}\gamma_{j,i}C^\omega$ for all $j\in J$, where all the $\gamma_{j,i}C^\omega$ may be assumed pairwise disjoint. So $W=\cup_{i\in I,j\in J}W_j\cap\gamma_{j,i}C^\omega$ is an open-separated union witnessed by $W_j\cap\gamma_{j,i}C^\omega\subseteq\gamma_{j,i}C^\omega$. Since all the subclasses of the quasi-Borel sets are closed under intersection with clopen sets by Lemma 1.1 in~\cite{Martin90}, the $W_j\cap\gamma_{j,i}C^\omega$ are all below $W$ in the hierarchy, so each $\iota(W_j\cap\gamma_{j,i}C^\omega)$ is quasi-Borel by induction hypothesis. Therefore the open-separated union $\iota(W)=\cup_{i\in I,j\in J}\iota(W_j\cap\gamma_{j,i}C^\omega)$ that is witnessed by $\iota(W_j\cap\gamma_{j,i}C^\omega)\subseteq\iota(\gamma_{j,i}C^\omega)\subseteq\iota(\gamma_{j,i})C^\omega$ (by property of $\iota$) is also quasi-Borel.\qed\smallskip
\end{iteMize}

\noindent Now let us justify explicitly why Theorem~\ref{thm:borel-nash} is, indeed, a generalisation of quasi-Borel determinacy. Let $T$ be a tree and $G(W;T)$ be a game as in~\cite{Martin90}, where $W$ is a quasi-Borel set of $[T]$, the infinite paths of the tree. Let $A:=\{a,b\}$ be a two-element set; let $C$ be the set of all elements that occur in the sequences of $T$, let $d$ map $C^{2*}$ to $a$ and $C^{2*+1}$ to $b$; let $O:=\{w,l\}$ be a two-element set; let $v$ map to $w$ all the elements of $W$ union "fatal" clopen sets that make $b$ lose when playing outside of $T$, and let $v$ map the rest of $C^\omega$ to $l$; finally, let $\prec_a:=\{(l,w)\}$ and  $\prec_b:=\{(w,l)\}$. The preferences are obviously strictly well-founded, $v^{-1}(\{w\})$ is quasi-Borel, so the game $\langle A,C,d,O,v,(\prec_a)_{a\in A}\rangle$ has a Nash equilibrium by Theorem~\ref{thm:borel-nash}. Restricting this Nash equilibrium to $T$ yields a strategy profile for $G(W;T)$, and one of the strategies involved is a winning strategy. 

Furthermore, since the notions of strict well-foundedness and acyclicity coincide on finite domains, Theorem~\ref{thm:borel-nash} is also a proper generalisation of the result~\cite{SLR-PhD08}\cite{SLR09} already mentioned in the abstract and in the end of Section~\ref{subsect:pw}.

Eventually, the second application of Theorem~\ref{thm:wo-ne} is the following corollary, assuming that the axiom of determinacy (AD) holds for every countable $C$. The proof of $1\Rightarrow 2$ may invoke Lemma~\ref{lem:no-gsg} and the proof of $2\Rightarrow 1$ may be copied from the proof of Theorem~\ref{thm:borel-nash}.

\begin{thm}\label{thm:det-nash}
Assume AD. Let $A$ and $O$ be non-empty sets with $O$ countable, and let $C$ be a countable set containing at least two elements. For every $a\in A$ let $\prec_a$ be a binary relation over $O$. 
Then the following two propositions are equivalent. 
\begin{enumerate}[\em(1)]
\item $\prec_a^{-1}$ is strictly well-founded for all $a\in A$.
\item Every infinite sequential game $\langle A,C,d,O,v,(\prec_a)_{a\in A}\rangle$ has a Nash equilibrium. 
\end{enumerate}
\end{thm}

\section*{Acknowledgement}

I thank Achim Blumensath, Benno van den Berg, and Yiannis N. Moschovakis for discussions on Borel determinacy. I also thank Arno Pauly for helpful discussions, comments, and ideas on drafts of the paper. Finally, I am especially grateful to Vassilios Gregoriades for his very helpful set-theoretic explanations and suggestions and to an anonymous referee for very detailed comments ranging from tiny typo corrections to insightful suggestions.

\section{Conclusion}

Theorem~\ref{thm:borel-nash} generalises quasi-Borel determinacy in a game-theoretic direction: it presents a general setting where strict well-foundedness is a necessary and sufficient condition for existence of Nash equilibrium. This should be of interest to game theorists.

Theorem~\ref{thm:borel-nash} does not generalise quasi-Borel determinacy in a set-theoretic direction, though, since it involves only quasi-Borel sets. However, Theorem~\ref{thm:wo-ne} says that the existence of a Nash equilibrium in a game involving sets of a given complexity is a consequence of determinacy for these sets. Said otherwise, Theorem~\ref{thm:wo-ne} suggests that using the concept of Nash equilibrium, instead of determinacy, is not likely to help refine the description of the complexity of sets, at least under well-foundedness assumptions. In this respect, Theorem~\ref{thm:wo-ne} may be seen a negative result by set-theorists.

There may be a positive set-theoretic side to this, nonetheless: since (quasi-)Borel determinacy was applied in several areas, it would be interesting to see how these applications may be generalised by invoking, \textit{e.g.}, Theorem~\ref{thm:borel-nash} instead of (quasi-)Borel determinacy.

\bibliographystyle{plain}     
\bibliography{article}

\begin{thebibliography}{10}

\bibitem{Borel21}
Emile Borel.
\newblock La th\'eorie du jeu et les equations int\'egrales \`a noyau
  sym\'etrique.
\newblock {\em Comptes Rendus Hebdomadaires des s\'eances de l'acad\'emie des
  Sciences}, 173:1304--1308, 1921.

\bibitem{GS53}
D.~Gale and F.~M. Stewart.
\newblock Infinite games with perfect information.
\newblock {\em Annals of Math. Studies}, 28:245--266, 1953.

\bibitem{Kuhn53}
Harold~W. Kuhn.
\newblock Extensive games and the problem of information.
\newblock {\em Contributions to the Theory of Games II}, 1953.

\bibitem{SLR-PhD08}
St{\'e}phane Le~Roux.
\newblock {\em Generalisation and formalisation in game theory}.
\newblock Ph.{D}. thesis, Ecole Normale Supérieure de Lyon, January 2008.

\bibitem{SLR09}
St{\'e}phane Le~Roux.
\newblock Acyclic preferences and existence of sequential {N}ash equilibria: a
  formal and constructive equivalence.
\newblock In {\em TPHOLs, International Conference on Theorem Proving in Higher
  Order Logics}, Lecture Notes in Computer Science, pages 293--309. Springer,
  August 2009.

\bibitem{Martin75}
Donald~A. Martin.
\newblock Borel determinacy.
\newblock {\em Annals of Mathematics}, 102:363–371, 1975.

\bibitem{Martin90}
Donald~A. Martin.
\newblock An extension of {B}orel determinacy.
\newblock {\em Annals of Pure and Applied Logic}, 49:279--293, 1990.

\bibitem{Moschovakis09}
Yiannis~N. Moschovakis.
\newblock {\em Descriptive Set Theory}.
\newblock American Mathematical Society, 2009.

\bibitem{NM44}
John~von Neumann and Oskar Morgenstern.
\newblock {\em Theory of Games and Economic Behavior}.
\newblock Princeton Univ. Press, Princeton, 1944.

\bibitem{OR94}
Martin~J. Osborne and Ariel Rubinstein.
\newblock {\em A Course in Game Theory}.
\newblock The MIT Press, 1994.

\bibitem{PS09}
Soumya Paul and Sunil Simon.
\newblock Nash equilibrium in generalised muller games.
\newblock In Ravi Kannan and K~Narayan Kumar, editors, {\em IARCS Annual
  Conference on Foundations of Software Technology and Theoretical Computer
  Science (FSTTCS 2009)}, volume~4 of {\em Leibniz International Proceedings in
  Informatics (LIPIcs)}, pages 335--346, Dagstuhl, Germany, 2009. Schloss
  Dagstuhl--Leibniz-Zentrum fuer Informatik.

\bibitem{Ummels05}
Michael Ummels.
\newblock Rational behaviour and strategy construction in infinite multiplayer
  games, 2005.
\newblock Master's thesis supervised by Erich Gr\"adel.

\bibitem{Wolfe55}
Philip Wolfe.
\newblock The strict determinateness of certain infinite games.
\newblock {\em Pacific Journal of Mathematics}, 5:841–847, 1955.

\end{thebibliography}

\end{document}